\documentclass[journal]{IEEEtran}

\usepackage{amssymb}
\usepackage{amsmath}
\usepackage{amsthm}
\usepackage{array}
\usepackage{xcolor}
\usepackage{float}
\usepackage[caption=false,font=footnotesize]{subfig}
\usepackage{graphicx}
\usepackage{tikz}
\usetikzlibrary{positioning, arrows.meta, calc, fit}
\usepackage{cite}

\floatstyle{ruled}
\newfloat{algorithm}{tbp}{loa}
\providecommand{\algorithmname}{Algorithm}
\floatname{algorithm}{\protect\algorithmname}

\newtheorem{example}{Example}
\newtheorem{definition}{Definition}
\newtheorem{prob}{Problem}
\newtheorem{lem}{Lemma}
\newtheorem{prop}{Proposition}
\newtheorem{thm}{Theorem}
\newtheorem{cor}{Corollary}
\newtheorem{remark}{Remark}

\newcommand{\V}{\mathcal{V}}
\newcommand{\E}{\mathcal{E}}
\newcommand{\G}{\mathcal{G}}
\newcommand{\U}{\mathcal{U}}
\newcommand{\Y}{\mathcal{Y}}
\newcommand{\Hinf}{\mathcal{H}_{\infty}}
\newcommand{\R}{\mathbb{R}}
\newcommand{\diag}{\mathrm{diag}}
\begin{document}

\title{Interface-Preserving $\mathcal{H}_\infty$\\ Model Reduction of Network Systems}

\author{Noam~Leiter and Daniel~Zelazo~\IEEEmembership{Senior Member,~IEEE}%
\thanks{N. Leiter is with Prime Minute, Palo Alto, California 94301,
        (e-mail: noam.leiter@gmail.com).}%
\thanks{D. Zelazo is with the Stephen B. Klein Faculty of Aerospace Engineering,
        Technion--Israel Institute of Technology, Haifa 32000, Israel
        (e-mail: dzelazo@technion.ac.il). \emph{Corresponding author.}}%
\thanks{This research was supported by Grant No.~2285/20 from the Israel
        Science Foundation.}%
}

\maketitle

\begin{abstract}
We study projection-based reduced-order models that exactly retain specified
input and output interfaces.  Admissible projections are characterized
geometrically: the retained subspace must contain every input and output
direction, which fixes the minimum order.  The reduction error factors into
input-coupling, internal-dynamics, and output-coupling systems, giving a
submultiplicative $\mathcal{H}_\infty$ bound; for positively dominated factors
the $\mathcal{H}_\infty$ error equals a dc gain computable from linear solves.
Strict dissipativity certifies the stability of every orthogonal reduction,
and a Hankel singular value bounds the optimum a posteriori.  These results
motivate sequential graph-coarsening heuristics for networks.  On a
$65{,}566$-state chip thermal network, interface-preserving coarsening yields
sparse, positive lumped models that are more accurate than Kron reduction of
the same interface at every order, with far fewer edges; balanced truncation
is more accurate but destroys the network structure.
\end{abstract}


\IEEEpeerreviewmaketitle

\section{Introduction}\label{sec.intro}

Model reduction provides low-order surrogates for the analysis and control of large-scale systems \cite{dullerud2000course}.  Balanced truncation yields a stable reduced model with an a priori $\mathcal{H}_\infty$ error bound \cite{moore1981principal,antoulas2005approximation} and is routine at $10^4$--$10^5$ states with low-rank Lyapunov solvers \cite{simoncini2016computational}, while Krylov and other projection methods scale further but do not guarantee stability for an arbitrary projection \cite{antoulas2005approximation}.  For network systems, however, the drawback of these methods is structural: the reduced coordinates mix all nodes, so the reduced model loses its graph interpretation and its inputs and outputs no longer correspond to physical nodes.

Structured reductions of networks preserve the graph topology, which encodes stability, controllability, and synchronization properties.  Graph partitions have been used to reduce consensus, passive, and dissipative networks \cite{ishizaki2014model,monshizadeh2014projection,besselink2016clustering}, with clusters selected, for example, by the dissimilarity of nodal responses \cite{cheng2019dissimilarity}; see \cite{cheng2021networkreview} for a survey.  Clustering of positive directed networks that preserves positivity and steady-state behavior \cite{ishizaki2015clustered} is the closest prior work to the positively dominated setting considered here.  Graph-based $\mathcal{H}_\infty$ reduction of multi-agent systems is studied in \cite{leiter2017graph}, a priori $\mathcal{H}_2$ and $\mathcal{H}_\infty$ error bounds for clustering of leader--follower systems are derived in \cite{jongsma2018model}, and structure-preserving balanced truncation of networked passive systems is developed in \cite{cheng2019balanced}.  More recently, the edge weights of a clustered consensus network have been optimized for $\mathcal{H}_2$ or $\mathcal{H}_\infty$ error \cite{sabbagh2022consensus}, structure-preserving reduction has been extended to interconnected systems \cite{vandewouw2025interconnected}, and an $\mathcal{H}_2$-optimal Petrov--Galerkin framework preserves stability and dissipativity for every projection basis by weighting the projection with a Lyapunov certificate, without requiring specified input and output nodes to be retained \cite{cheng2026optimal}.

In this work we require every actuated and measured node to retain its individual identity in the reduced model, a stronger condition than the approximate interface matching permitted by most clustering methods.  This is the projection-based counterpart of Kron reduction \cite{dorfler2013kron}, which keeps a set of boundary nodes and eliminates the interior through a Schur complement.  Kron reduction is exact at dc but carries no error guarantee away from dc, and it densifies the reduced graph, since boundary nodes joined through the interior become adjacent.  An interface-invariant reduced model instead aggregates the interior into a coarse graph, so interior dynamics remain in averaged form and the error is controlled in $\mathcal{H}_\infty$.  Within this family we ask which reduced orders are achievable, how the $\mathcal{H}_\infty$ error can be bounded, and how the remaining projection freedom should be used.  The central tool is a product factorization of the error system that separates input coupling, discarded dynamics, and output coupling.  For positively dominated systems it yields an exact expression: the $\mathcal{H}_\infty$ error equals a dc gain that can be evaluated without frequency gridding, even at tens of thousands of states.

The main contributions are as follows.  (i)~We characterize all interface-invariant orthogonal projections geometrically: the retained subspace must span every input and output direction, which establishes the minimum reduced order and the remaining degrees of freedom.  (ii)~We derive a product representation of the projection-based reduced-order model (PROM) $\mathcal{H}_\infty$ error, with explicit stability hypotheses; for interface-invariant PROMs the coupling factors are strictly proper.  (iii)~Although the reduction error system is not positive in general, we show that for positively dominated systems \cite{rantzer2012optimizing} its $\mathcal{H}_\infty$ norm reduces to a closed-form steady-state gain, as for internally positive systems \cite{tanaka2011bounded}; for strictly dissipative systems, a projection-dependent upper bound follows from linear solves, without a frequency sweep.  (iv)~We express a projection as a sequence of singleton projections, characterize when the resulting triangle-inequality bound is tight, and use it to motivate sequential graph-coarsening heuristics, with the $(r+1)$st Hankel singular value as an a posteriori approximation factor.  (v)~On a $65{,}566$-state processor thermal network with $30$ functional-unit interfaces, we compare interface-preserving graph coarsening with Kron reduction, coordinate deletion, and balanced truncation in accuracy, preserved structure, and cost.

Section \ref{sec.problem} formulates the reduction problem and interface invariance.  Section \ref{sec.structure} gives the singleton sequence and product factorization.  Section \ref{sec.positive} develops the positively dominated specialization.  Section \ref{sec.MAS_IIPROM} presents structured sequential heuristics, Section \ref{sec.cases} presents the chip thermal case study, and Section \ref{sec.conclusion} concludes.

\paragraph*{Notations}

We write $A\succ0$ for positive-definite, $A\geq0$ for entrywise nonnegative, and call $A$ \emph{Metzler} if $[A]_{ij}\geq0$ for all $i\neq j$. For symmetric $A\in\R^{n\times n}$, $\lambda(A)=\{\lambda_k(A)\}_{k=1}^n$ denotes the eigenvalues in ascending order. $[A]_{ij}$ is the $(i,j)$ entry, $\diag(A,B)$ a block-diagonal matrix, and $\|A\|_2$ the spectral norm. The $i$th standard basis vector of $\R^n$ is $e_i^n$. For an index set $\mathcal{I}=\{i_1,\ldots,i_m\}$, $e_{\mathcal{I}}^n=[e_{i_1}^n\ \cdots\ e_{i_m}^n]$. $\operatorname{im}M$ is the column space of $M$. $\mathcal{RH}_\infty^{p\times m}$ denotes the set of stable proper real-rational $p\times m$ transfer matrices. A graph $\mathcal{G}=(\mathcal{V},\mathcal{E})$ has vertex set $\mathcal{V}$ and edge set $\mathcal{E}$. Nodes $u,v$ are \emph{adjacent} ($u\sim v$) if $uv\in\mathcal{E}$.

\section{Model Reduction and Problem Formulation}\label{sec.problem}

In this section we review some fundamental notions from the theory of model reduction.  We then develop a specialization of the model reduction problem for multi-agent systems.  With this framework we will then formally state the problem we aim to solve.

\subsection{The Model Reduction Problem}\label{subsec.modelred}

We consider an LTI system $\Sigma$ with minimal state-space realization
\begin{align}\nonumber
\Sigma \,:\, \begin{cases} \dot x(t) &=Ax(t) + Bu(t) \\ y(t) &= Cx(t) + Du(t), \end{cases}
\end{align}
where $x(t) \in \R^{n_x}$ is the system state, $u(t) \in \R^{n_u}$ is the input, and $y(t) \in \R^{n_y}$ is the output.  We denote the state-space realization of a system by the tuple $(A,B,C,D)$.  The transfer function matrix (TFM) of the system $\Sigma$ is given as
\begin{align}\nonumber
\hat \Sigma(s) = C(sI_{n_x}-A)^{-1}B+D.
\end{align}
Hereafter, the notation $\hat \Sigma$, without explicit dependence on $s$, will be used to denote the TFM of a system $\Sigma$.
A realization $\Sigma:=(A,B,C,D)$ is \emph{minimal} if it is controllable and observable.  The \emph{order} of a system is its \emph{McMillan degree}, denoted as $\deg(\Sigma)$.

The model reduction problem is to find an order-$r$ system $\Sigma_r$, with $r<n_x$, whose transfer matrix approximates $\hat\Sigma$.  The reduced and full models retain the same external input and output dimensions: $B_r\in\R^{r\times n_u}$ and $C_r\in\R^{n_y\times r}$.  We define the error system
\begin{align}\label{error_system}
\Sigma_e = \Sigma - \Sigma_r,
\end{align}
which has realization
\begin{align}\label{error_system_ss}
    \Sigma_e \,:\, &= (A_e,B_e,C_e,D_e) \nonumber\\
    &= \left(\begin{bmatrix} A & 0 \\ 0 & A_r\end{bmatrix}, \begin{bmatrix} B \\ B_r\end{bmatrix},\begin{bmatrix} C & -C_r \end{bmatrix},D-D_r \right).
\end{align}
Thus, one can consider, for a fixed $r$, finding a system $\Sigma_r$ such that $\|\Sigma_e\|$ is minimized.  For a detailed review see \cite{antoulas2005approximation}. In this work we focus on $\Hinf$ reduction errors.  For a stable and proper system,
$$\|\hat \Sigma\|_\infty = \sup_{\omega \in \R} \bar \sigma(\hat \Sigma(j\omega)),$$
where $\bar \sigma(M)$ is the largest singular value of the matrix $M$.

A widely used family of reduction methods are the so-called \emph{projection-based reductions}.  A Petrov--Galerkin reduction uses trial and test bases $V,W\in\R^{n_x\times r}$ satisfying $W^TV=I_r$ and gives
\begin{align}\label{PROM}
\Sigma_r(W,V) :=(W^TAV,W^TB,CV,D).
\end{align}
When $W=V=P$ and $P^TP=I_r$, the reduction is orthogonal.  We denote the resulting realization by $\Sigma_r(P)$; it is obtained from \eqref{PROM} by setting $W=V=P$.  We call this model an \emph{orthogonal projected reduced-order model} (PROM).

\subsection{Multi-Agent Systems and Interface-Invariant Model Reduction}\label{subsec.masred}
\begin{figure}[!t]
    \centering
    \includegraphics[width=0.4\textwidth]{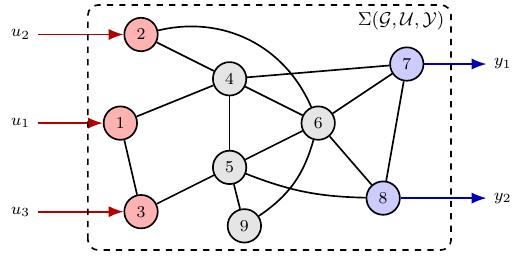}
    \caption{An MAS with control nodes ({\textcolor{red}{red}}) and output nodes ({\textcolor{blue}{blue}}).}\label{fig.MAS}
\end{figure}

We consider linear multi-agent systems (MAS) comprised of $n$ agents interacting over a network $\G=(\V,\E)$, with state $x_i\in\R^{q_x}$ per agent ($nq_x=n_x$). The block structure of $A(\G)$ follows the graph: the $ij$-block is nonzero iff $ij\in\E$. Input nodes $\U\subseteq\V$ ($|\U|=m$) receive external inputs $u_i\in\R^{q_u}$; output nodes $\Y\subseteq\V$ ($|\Y|=p$) are measured, $y_i\in\R^{q_y}$. The MAS realization is
\begin{align}\label{MAS_realization}
\Sigma(\G,\U,\Y) := (A(\G),B(\U),C(\Y),D(\U,\Y)),
\end{align}
with $B(\U)\in\R^{n_x\times n_u}$, $n_u=mq_u$, and $C(\Y)\in\R^{n_y\times n_x}$, $n_y=pq_y$ (Figure~\ref{fig.MAS}). Agents in $\U\cup\Y$ are the \emph{interface agents}; self-loops are permitted.

In this work we will examine a special class of PROMs that we term \emph{interface-invariant PROMs}. These are reduced models that maintain the  input-output  structure  of  the  full  order  model under the projection operation. Such PROMs are required, for example, for the reduction of the MAS \eqref{MAS_realization} where the interface agent structure is maintained in the reduced model.
We now define formally the notion of an interface-invariant PROM (IIPROM).
\begin{definition}[IIPROM]\label{IIPROM} Given a system with realization
$\Sigma:=\left(A,B,C,D\right)$, an \emph{interface invariant projection-based reduced order model} (IIPROM) of $\Sigma$ is any PROM $\Sigma_{r}(P)$
such that $C=CPP^T$ and $B=PP^T B$.
\end{definition}
Note that an IIPROM does not require $PP^T=I_{n_x}$. For example, if $C=B^T=[1\ 0\ 0]$, the orthonormal matrix
$$P=\begin{bmatrix}1&0&0\\0&1/\sqrt2&1/\sqrt2\end{bmatrix}^T$$
satisfies $C=CPP^T$ and $B=PP^TB$.

\begin{remark}
Although motivated by multi-agent systems, the IIPROM framework is not intrinsically graph based.  It applies to any orthogonal PROM for which prescribed input and output state-space directions must be retained, regardless of whether the system has a network interpretation.
\end{remark}

Define the \emph{protected input--output subspace}
\begin{equation}
\mathcal S_{\rm io}:=\operatorname{im}B+\operatorname{im}C^T
=\operatorname{im}\begin{bmatrix}B&C^T\end{bmatrix},
\;\; k:=\dim\mathcal S_{\rm io}.\nonumber
\end{equation}
The following result separates the interface requirement from the remaining model-reduction freedom.

\begin{thm}[Geometric characterization of IIPROMs]\label{thm:IIPROM-geometry}
Let $P\in\R^{n_x\times r}$ satisfy $P^TP=I_r$, and let
$\Pi=P P^T$.  Then $\Sigma_r(P)$ is interface invariant if and only if
\begin{equation}
\mathcal S_{\rm io}\subseteq\operatorname{im}P.
\label{eq:interface-subspace-containment}
\end{equation}
Consequently, an order-$r$ interface-invariant projection exists if and
only if $k=\dim\mathcal S_{\rm io}\leq r\leq n_x$, and it is a strict reduction precisely when
$k\leq r<n_x$.

More precisely, let $U_{\rm io}\in\R^{n_x\times k}$ be an
orthonormal basis for $\mathcal S_{\rm io}$, and complete it to the
orthogonal matrix $[\,U_{\rm io}\ U_\perp\,]$.  Every admissible
retained subspace has an orthonormal basis $P$ of the form
\begin{equation}
P=\begin{bmatrix}U_{\rm io}&U_\perp Z\end{bmatrix}R,
\,
Z^TZ=I_{r-k},\, R^TR=RR^T=I_r,
\label{eq:IIPROM-parameterization}
\end{equation}
where $Z\in\R^{(n_x-k)\times(r-k)}$ and $R\in\R^{r\times r}$ is an orthogonal basis-change within the retained subspace; since $\Pi=PP^T$ is invariant under right-multiplication by $R$, the factor $R$ does not affect the reduced model.
Conversely, every matrix in \eqref{eq:IIPROM-parameterization} defines
an order-$r$ IIPROM.
\end{thm}

\begin{proof}
Because $\Pi$ is the orthogonal projector onto $\operatorname{im}P$,
$\Pi B=B$ holds if and only if $\operatorname{im}B\subseteq
 \operatorname{im}P$ (similarly condition for $\Pi C^T=C^T$).
The two interface equalities are therefore equivalent to
\eqref{eq:interface-subspace-containment}, which also proves the order
condition.

For the parameterization, suppose \eqref{eq:interface-subspace-containment}
holds.  The retained subspace has the orthogonal decomposition
\[
\operatorname{im}P=\mathcal S_{\rm io}\oplus\mathcal W,
\qquad
\mathcal W\subseteq\mathcal S_{\rm io}^{\perp},\qquad
\dim\mathcal W=r-k.
\]
Writing an orthonormal basis of $\mathcal W$ in the coordinates of
$U_\perp$ gives $U_\perp Z$ with $Z^TZ=I_{r-k}$.  Thus
$P_0=[\,U_{\rm io}\ U_\perp Z\,]$ and $P$ are orthonormal bases for
the same subspace, so $P=P_0R$ for an orthogonal matrix $R$.
The converse follows because every such $P$ contains
$\mathcal S_{\rm io}$ in its range.
\end{proof}

Right multiplication by $R$ only changes reduced coordinates.
The matrices $\Sigma_r(P_0R)$ and $\Sigma_r(P_0)$ are similar and have the same
transfer matrix.  Hence one may set $R=I_r$ without loss of input--output
generality.  The constrained search over $P$ is thereby reduced to the
choice of an $(r-k)$-dimensional subspace of
$\mathcal S_{\rm io}^{\perp}$, represented by $Z$ in
\eqref{eq:IIPROM-parameterization}.  For a multi-agent system with
coordinate-selecting $B$ and $C$, $\mathcal S_{\rm io}$ is simply the
span of the actuated and measured agent coordinates.

With the above notions in place, we can now formally state the optimal IIPROM problem.

\begin{prob}[optimal IIPROM]\label{prob:optimalIIPROM}
Consider a stable proper system of order $n_x$ with realization $\Sigma:=\left(A,B,C,D\right)$.
 Find $P\in\R^{n_x\times r}$ with $P^T P=I_{r}$ such that
the PROM $\Sigma_{r}(P)$ is interface invariant and stable and
minimizes the $\Hinf$-norm of the reduction error system \eqref{error_system},
\begin{align}
\min_{P\in\mathbb{R}^{n_x\times r}} & \| \Sigma_{e}\| _{\infty}\nonumber\\
\mathrm{s.t.} \quad & P^T P=I_{r},\quad C=CPP^T,\nonumber\\
                    & B=PP^T B,\quad P^TAP\ \text{Hurwitz}.\nonumber
\end{align}
\end{prob}

Theorem \ref{thm:IIPROM-geometry} removes the interface-invariance constraints from Problem~\ref{prob:optimalIIPROM} by replacing $P$ with the parameterization \eqref{eq:IIPROM-parameterization}, reducing the search to the choice of $Z\in\R^{(n_x-k)\times(r-k)}$ with orthonormal columns.  The remaining $\mathcal H_\infty$ objective is nevertheless nonsmooth and nonconvex, and reduced-model stability remains a nonlinear constraint.  Consequently, the problem generally calls for optimization on the reduced Grassmann manifold associated with $\mathcal S_{\rm io}^{\perp}$ or for a structured restriction of the feasible projections; a Riemannian method of this kind for structure-preserving $\mathcal{H}_2$ reduction is developed in \cite{cheng2026optimal}.

While a closed-form solution is not available in general, there are instances where an analytic solution exists. We provide a simple example below.

\begin{figure}[!t]
\centering
\includegraphics[width=0.88\linewidth]{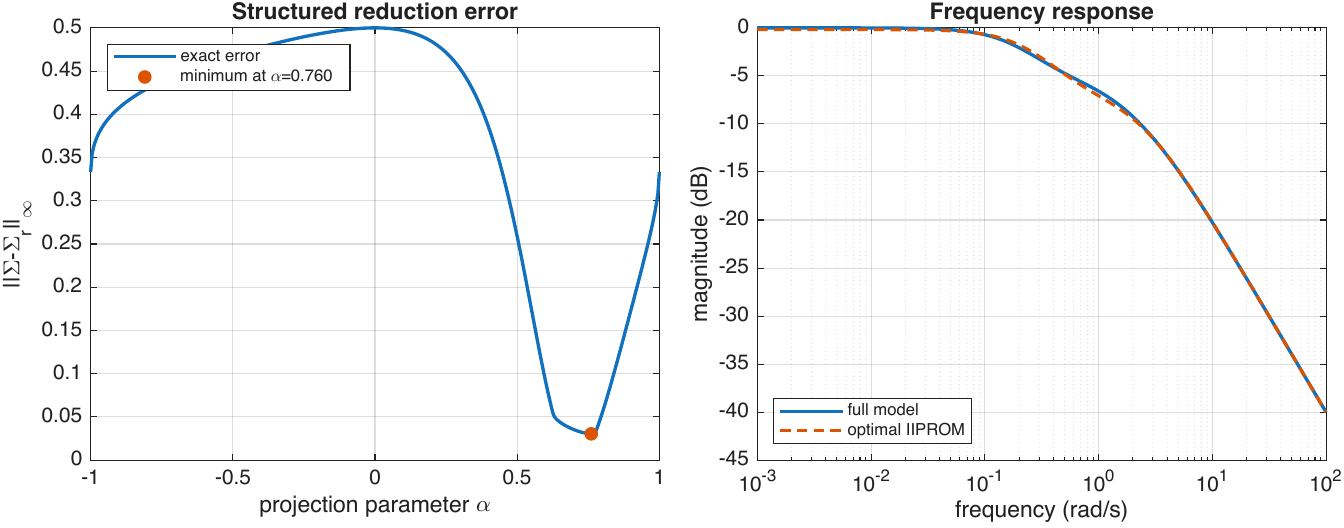}
\caption{Example~\ref{IIPROMexample}: reduction error $\|\Sigma_e\|_\infty$ vs.\ $\alpha$ with its minimizer (left), and the full and optimal reduced frequency responses (right).}
\label{fig:introductory-example}
\end{figure}

\begin{example}\label{IIPROMexample}
Consider the SISO system $\Sigma=(A,B,C,D)$ with $A=\left[\begin{smallmatrix}-2&1&0\\1&-2&1\\0&1&-1\end{smallmatrix}\right]$, $B=C^T=e_1$, and $D=0$, whose transfer function is $\hat\Sigma=(s^2+3s+1)/(s^3+5s^2+6s+1)$, with $\|\Sigma\|_\infty=1$.  We seek an IIPROM of order $r=2$.  Every $P\in\mathbb{R}^{3\times2}$ with $P^TP=I_2$, $CPP^T=C$, and $PP^TB=B$ has the form $P(\alpha)=[e_1\ \ \alpha e_2+\beta(\alpha)e_3]$, with $\alpha\in[-1,1]$ and $\beta(\alpha)=\sqrt{1-\alpha^2}$.  The resulting IIPROM is
\[
\hat\Sigma_r=\frac{s+\gamma^2}{s^2+(\gamma^2+2)s+2\gamma^2-\alpha^2},
\]
with $\gamma^2=\alpha^2+(\alpha-\beta(\alpha))^2$.  It is stable for all $\alpha$ and minimal for $\alpha\neq0$; at $\alpha=0$ the retained direction $e_3$ is decoupled from $e_1$ and $\hat\Sigma_r=1/(s+2)$.  The error $\|\Sigma_e\|_\infty$ ranges from $0.5$ at $\alpha=0$ to its minimum $\approx0.03$ at $\alpha^*=0.76$, which solves Problem~\ref{prob:optimalIIPROM}, and it is nonconvex and nonsmooth in $\alpha$ (Fig.~\ref{fig:introductory-example}).
\end{example}

Example \ref{IIPROMexample} exploits the fact that the free subspace in Theorem \ref{thm:IIPROM-geometry} is one dimensional.  In general the geometric parameterization removes the two interface equalities but does not produce a closed-form or convex solution of Problem \ref{prob:optimalIIPROM}: the choice of the complementary subspace, the $\mathcal H_\infty$ objective, and the reduced stability condition remain coupled.  In the sequel, we investigate the IIPROM error structure and derive tractable suboptimal constructions.

\section{Structure and Analysis\\ of PROMs and IIPROMs}\label{sec.structure}

Understanding the structure of IIPROMs is essential for tackling Problem \ref{prob:optimalIIPROM} for multi-agent systems.  In this direction, we provide complementary perspectives on the error system associated with IIPROMs. First, we introduce a sequential decomposition of the projection matrix used in PROMs into a product of elementary “singleton” projections, each reducing the system order by one. This construction leads to a telescoping sum representation of the error system and yields a tractable upper bound on the overall reduction error. Building on this result, we formulate a suboptimal IIPROM design problem aimed at minimizing this bound. These results provide a foundation for algorithmic strategies that sequentially eliminate state variables while preserving input-output invariance and stability.

\subsection{The PROM Sequence Bound}

Stability of the full-order model alone does not guarantee stability of an orthogonal PROM.  The following proposition gives two sufficient conditions under which an orthogonal compression is stable. It will be used both for the singleton sequence below and later for graph-contraction candidates.

\begin{prop}[A priori stability of orthogonal compressions]\label{prop:compression-stability}
Let $P^TP=I_r$, $\Pi=PP^T$, and $A_P=P^TAP$.
\begin{enumerate}
\item If, for some $\mu>0$,
\begin{equation}
\frac{A+A^T}{2}\preceq-\mu I_n,
\label{eq:euclidean-dissipativity}
\end{equation}
then $A_P$ is Hurwitz for every $P$, with
$\max\operatorname{Re}\lambda(A_P)\leq-\mu$.
\item More generally, suppose $X=X^T\succ0$ satisfies
$A^TX+XA\prec0$.  If $X\Pi=\Pi X$, then $A_P$ is Hurwitz with Lyapunov matrix $X_P=P^TXP\succ0$.
\end{enumerate}
\end{prop}
\begin{proof}
For the first claim,
\begin{equation*}
\frac{A_P+A_P^T}{2}
=P^T\frac{A+A^T}{2}P\preceq-\mu I_r,
\end{equation*}
which bounds the real part of every eigenvalue of $A_P$.  For the second claim, commutation gives $XP=PX_P$ and $P^TX\Pi=P^TX$.  Hence
\begin{equation*}
A_P^TX_P+X_PA_P=P^T(A^TX+XA)P\prec0.
\end{equation*}
\end{proof}

In particular, under part~(i), every intermediate compression in any singleton projection sequence is Hurwitz with the same decay margin $\mu$.  More generally, part~(ii) applies whenever the chosen retained subspace is compatible with a Lyapunov matrix through the commutation condition.  Some additional assumption beyond mere Hurwitz stability of $A$ is essential, since a Hurwitz nonsymmetric matrix can have an unstable orthogonal compression.  The commutation condition can be avoided by replacing the orthogonal projection with the $X$-orthogonal Petrov--Galerkin projection $P(P^TXP)^{-1}P^TX$, which preserves stability, and more generally any dissipativity certified by $X$, for every $P$ \cite{cheng2026optimal}.  We retain orthogonal projections because the product form of Section~\ref{sec:The-Product-Form} and the quotient-graph interpretation of the partition projections in Section~\ref{sec.MAS_IIPROM} rely on them.

In the following section we present a lemma showing
that any projection from $\R^{n_x}$ to $\R^{r}$ can be obtained from
a sequence of $n_x-r$ projections, each reducing the dimension by one.
We denote such projections as \emph{singleton projections}.

\begin{lem}[Singleton factorization]\label{lem:singleton}
Let $P\in\R^{n_x\times r}$ with $P^T P=I_{r}$ and $r<n_x$. Then there exists a sequence $\{ P_{(k)}\}_{k=1}^{n_x-r}$ with
$P_{(k)}\in\R^{n_x-k+1\times n_x-k}$ and $P_{(k)}^T P_{(k)}=I_{n_x-k}$,
such that $P=\prod_{k=1}^{n_x-r}P_{(k)}$ (with factors ordered from left to right in increasing $k$).
\end{lem}
\begin{proof}
Let $P\in\R^{n_x\times r}$ with  $P^TP=I_{r}$. Then there exists $Q\in\R^{n_x\times n_x-r}$
with $Q^TQ=I_{n_x-r}$ such that $P^T Q=0$. Construct
\begin{align}\nonumber
 P_{(1)}=\begin{bmatrix}P&q_{1}&\cdots&q_{n_x-r-1}\end{bmatrix},
 \end{align}
where $q_{i}$ is the $i$th column of $Q$, and
\begin{align}\nonumber
P_{(k)}=\begin{bmatrix} I_{n_x-k} & 0_{n_x-k\times 1}\end{bmatrix}^T
\end{align}
for $k\in\left[2,\, n_x-r\right]$. Since $P^TP=I_{r}$ and $P^TQ=0$
 we have $P_{(1)}^TP_{(1)}=I_{n_x-1}$, and
 for $k\in\left[2,\, n_x-r\right]$ it follows that $P_{(k)}^TP_{(k)}=I_{n_x-k}$.
From this construction we get
\begin{align*}
\prod_{k=1}^{n_x-r}P_{(k)} &=P_{(1)}\left(\prod_{i=1}^{n_x-r-1}P_{\left(1+i\right)}\right)\\
 & =P_{(1)}\left(\prod_{i=1}^{n_x-r-1}\begin{bmatrix}I_{n_x-i}& 0_{(n_x-i)\times1}\end{bmatrix}^T\right) =P.
\end{align*}
\end{proof}

Note that one can generate other sequences $P=\prod_{k=1}^{n_x-r}\tilde{P}_{(k)}$
using the transformations $\tilde{P}_{(1)}=P_{(1)}U_{(1)}$,
$\tilde{P}_{(k)}=U_{\left(k-1\right)}^TP_{(k)}U_{(k)}$
for $k\in\left[2,\, n_x-r-1\right]$, and $\tilde{P}_{\left(n_x-r\right)}=U_{\left(n_x-r-1\right)}^TP_{\left(n_x-r\right)}$,
where $U_{(k)}\in\R^{(n_x-k)\times (n_x-k)}$ is an orthogonal
matrix.

By expressing a PROM as a sequence of singleton projections, we obtain
the following PROM sequence bound.

\begin{prop}[Singleton PROM sequence bound]\label{prop:seqbound}
Consider an LTI system $\Sigma_{(0)}:=\left(A_{(0)},B_{(0)},C_{(0)},D\right)$ and a PROM $\Sigma_{r}(P)$. Let $\left\{ P_{(k)}\right\} _{k=1}^{n_x-r}$
be a sequence with $P_{(k)}\in\R^{(n_x-k+1)\times (n_x-k)}$, $P_{(k)}^TP_{(k)}=I_{n_x-k}$
such that $P=\prod_{k=1}^{n_x-r}P_{(k)}$. Define recursively
\begin{align*}
A_{(k)}&=P_{(k)}^TA_{(k-1)}P_{(k)},& B_{(k)}&=P_{(k)}^TB_{(k-1)},\\
C_{(k)}&=C_{(k-1)}P_{(k)}.
\end{align*}
Assume every $A_{(k)}$, including $A_{(0)}$, is Hurwitz. Then the $\Hinf$
norm of the error reduction system, $\Sigma_{e}$ in \eqref{error_system},
is bounded by
\begin{equation}
\|\Sigma_{e}\|_\infty\leq\sum_{k=1}^{n_x-r}\|\Sigma_{(k)}-\Sigma_{\left(k-1\right)}\|_\infty,\label{eq:Singleton PROM sequence bound}
\end{equation}
where $\Sigma_{(k)}:=(A_{(k)},B_{(k)},C_{(k)},D)$.
\end{prop}
\begin{proof}
We express $\Sigma_{e}=\Sigma_{r}-\Sigma$
as the telescoping sum
$$\sum_{k=1}^{n_x-r}\left(\Sigma_{(k)}-\Sigma_{\left(k-1\right)}\right),$$
with $\Sigma_{(0)}=\Sigma$ and $\Sigma_{\left(n_x-r\right)}=\Sigma_{r}(P)$.
From the triangle inequality we obtain
$$\|\Sigma_{e}\|_\infty\leq\sum_{k=1}^{n_x-r}\|\Sigma_{(k)}-\Sigma_{\left(k-1\right)}\|_\infty,$$
completing the proof.
\end{proof}

The bound \eqref{eq:Singleton PROM sequence bound} can be conservative when the step errors peak at different frequencies or input directions.  The following result identifies the precise conditions under which equality holds.

\begin{prop}[Tightness of the singleton sequence bound]
\label{prop:sequence-tightness}
For $k=1,\ldots,n_x-r$, define the stable, strictly proper step error
\begin{equation}
E_k(s):=\hat\Sigma_{(k-1)}(s)-\hat\Sigma_{(k)}(s),
\qquad a_k:=\|E_k\|_\infty,\nonumber
\end{equation}
and omit any identically zero $E_k$.  Equality holds in
\eqref{eq:Singleton PROM sequence bound} if and only if there exist a
frequency $\omega_\star\in\R$ and unit vectors
$v_\star\in\mathbb C^{n_u}$ and $u_\star\in\mathbb C^{n_y}$ such that
\begin{equation}
E_k(j\omega_\star)v_\star=a_k u_\star
\quad\text{for every nonzero step }k.
\label{eq:sequence-alignment}
\end{equation}
Thus all steps must attain their $\mathcal H_\infty$ norms at one frequency,
on the same worst-case input direction, with phase-aligned output directions.
\end{prop}
\begin{proof}
The telescoping identity gives
$\hat\Sigma-\hat\Sigma_r=\sum_k E_k$.  If
\eqref{eq:sequence-alignment} holds, then
\begin{equation*}
\left\|\sum_kE_k(j\omega_\star)v_\star\right\|_2
=\sum_k a_k,
\end{equation*}
so the sequence bound is attained.  Conversely, the nonzero strictly proper
rational matrix $\sum_kE_k$ attains its $\mathcal H_\infty$ norm at a finite
frequency.  If the sequence bound is attained, choose a unit maximizing
right singular vector $v_\star$ at such a frequency $\omega_\star$.  Then
every inequality in
\begin{align*}
\left\|\sum_kE_k(j\omega_\star)v_\star\right\|_2
&\leq\sum_k\|E_k(j\omega_\star)v_\star\|_2\\
&\leq\sum_k\bar\sigma(E_k(j\omega_\star))
\leq\sum_k a_k
\end{align*}
must be an equality.  Equality in the Euclidean triangle inequality makes
the vectors $E_k(j\omega_\star)v_\star$ nonnegative multiples of a common
unit vector $u_\star$, while equality in the remaining inequalities makes
each multiple equal to $a_k$.  This is \eqref{eq:sequence-alignment}.
\end{proof}

When the step errors are positively dominated, the peak frequency of each $E_k$ is pinned to dc and the alignment condition of Proposition~\ref{prop:sequence-tightness} reduces to a purely algebraic criterion on the dc gains; this is explored in detail in Section~\ref{sec.positive}.

Proposition~\ref{prop:sequence-tightness} separates two sources of
conservatism in the sequence bound. Different peak frequencies prevent frequency
alignment, while different singular directions cause cancellation or dispersion
even when all steps peak at the same frequency.  The positively dominated
specialization developed in Section~\ref{sec.positive} removes the first source.

With Proposition \ref{prop:seqbound}, we can consider a suboptimal variation of Problem \ref{prob:optimalIIPROM} that attempts to minimize the upper bound for the reduction error.

\begin{prob}[suboptimal IIPROM sequence bound]\label{prob.sequencebound}
Consider a stable proper system $\Sigma$ of order $n_x$, with realization
$(A,B,C,D)$. Find a sequence $\{P_{(k)}\}_{k=1}^{n_x-r}$ such that the
associated PROM is interface invariant, all intermediate PROMs are stable,
and the sequence bound \eqref{eq:Singleton PROM sequence bound} is minimized:
\begin{align}
\min_{\{P_{(k)}\}} & \sum_{k=1}^{n_x-r}\|\Sigma_{(k)}-\Sigma_{\left(k-1\right)}\|_\infty\\
s.t.\,\,\,
 & P=\prod_{k=1}^{n_x-r}P_{(k)},\, P_{(k)}\in\R^{(n_x-k+1)\times (n_x-k)},\nonumber \\
 & P_{(k)}^TP_{(k)}=I_{n_x-k},\quad A_{(k)}\ \text{Hurwitz for all }k,\nonumber\\
 & C=CPP^T, \; B=PP^T B, \nonumber
\end{align}
where $\Sigma_{(k)}$ is defined recursively as in Proposition \ref{prop:seqbound}.
\end{prob}

Problem \ref{prob.sequencebound} is not a convex relaxation of Problem \ref{prob:optimalIIPROM}.  Its practical benefit is algorithmic: it permits a greedy search over admissible singleton reductions, with stability and interface invariance checked after every step.  The accumulated objective is a certified upper bound by the triangle inequality, but its conservatism can grow with the number of steps.

\subsection{The Product Form of Orthogonal PROMs }\label{sec:The-Product-Form}

It will also be useful to represent a PROM as a subsystem of $\Sigma$ after an appropriately defined transformation. We apply the unitary transformation $U=\begin{bmatrix}P & Q\end{bmatrix}$ on $\Sigma$, where $Q$ is the orthogonal complement of $P$, satisfying $Q^{T}P=0$ and $Q^{T}Q=I_{n_x-r}$. The resulting realization is
\begin{equation}
\tilde{\Sigma}:=\left(U^TAU,U^TB,CU,D\right),
\end{equation}
which can be expanded to
{\footnotesize{\begin{equation}
\tilde{\Sigma}:=\left(\left[\begin{array}{c|c}
A_{PP} & A_{PQ}\\
\hline A_{QP} & A_{QQ}
\end{array}\right],\left[\begin{array}{c}
P^{T}B\\
\hline Q^{T}B
\end{array}\right],\left[\begin{array}{c|c}
CP & CQ\end{array}\right],D\right),\label{eq:full order system PQ realization}
\end{equation}}}
where $A_{PP}\triangleq P^TAP$, $A_{PQ}\triangleq P^TAQ$, $A_{QP}\triangleq Q^TAP$ and $A_{QQ}\triangleq Q^TAQ$.

\begin{remark}[Interface invariance in PQ coordinates]\label{rem:iiprom-PQ}
If $\Sigma_{r}$ is an IIPROM, then $Q^{T}B=0$ and $CQ=0$ in \eqref{eq:full order system PQ realization}, since $Q^{T}B=Q^{T}(PP^{T}B)=0$ and $CQ=(CPP^{T})Q=0$. Consequently, the PQ realization of any IIPROM is block-lower-triangular: the retained-subspace dynamics are decoupled from the discarded inputs and outputs.
\end{remark}

The sequence bound of the previous subsection controls the cumulative error but treats each step error $\|\Sigma_{(k)}-\Sigma_{(k-1)}\|_\infty$ as an opaque quantity.  To sharpen the analysis, we need to understand the internal structure of the PROM error system. Specifically, how the input coupling, the discarded dynamics, and the output coupling each contribute.  We show that the error system of any orthogonal PROM factors as a cascade of three LTI systems that isolate these three contributions, reducing the full-order augmented realization of dimension $\deg(\Sigma)+\deg(\Sigma_r)$ to a product of systems of dimensions $\deg(\Sigma_r)$, $\deg(\Sigma)$, and $\deg(\Sigma_r)$.

\begin{thm}\label{thm:PROMbound}
Let $\Sigma:=(A,B,C,D)$, and consider a PROM $\Sigma_{r}(P)$. At every $s$ for which the inverses below exist, the error system $\Sigma_{e}=\Sigma-\Sigma_{r}(P)$ has the TFM
\begin{equation}
\hat{\Sigma}_{e}=C\Phi^{-1}Q(Q^T\Phi^{-1}Q)^{-1}Q^T\Phi^{-1}B,\label{eq:PROM error TFM}
\end{equation}
where $\Phi\triangleq sI_{n_x}-A$, and $Q$ satisfies $Q^T P=0$ and $Q^T Q=I_{n_x-r}$. Furthermore, $\Sigma_{e}$ can be expressed as the product of three systems,
\begin{equation}
\Sigma_{e}=\Theta(\Sigma,P)\Delta(\Sigma,P)\Gamma(\Sigma,P),\label{eq:PROM error realization}
\end{equation}
with realizations
\begin{subequations}\label{prod_form_elements}
\begin{align}
\Theta(\Sigma,P)&:=(A_{PP},A_{PQ},CP,CQ)\label{lti.1},\\
\Gamma(\Sigma,P)&:=(A_{PP},P^T B,A_{QP},Q^T B)\label{lti.2},\\
\Delta(\Sigma,P)&:=(A,Q,Q^T,0_{(n_x-r)\times(n_x-r)}).\label{lti.3}
\end{align}
\end{subequations}
\end{thm}
The proof of Theorem \ref{thm:PROMbound} is given in the Appendix. Because both sides are rational matrices, equality away from their finite pole sets establishes the transfer-matrix identity.
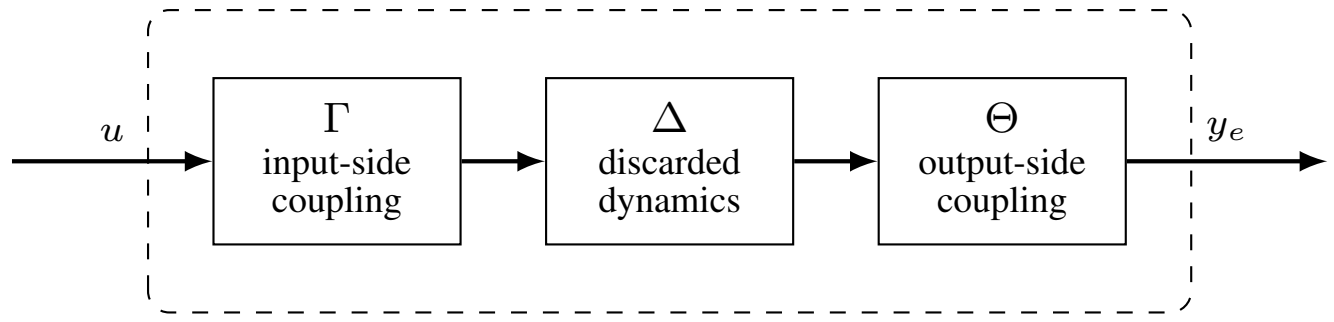
\begin{figure}[!t]
\centering
\begin{tikzpicture}[
    font=\footnotesize,
    >=latex,
    block/.style={draw, rectangle, minimum width=1.6cm, minimum height=1.1cm,
                  align=center, text width=1.4cm},
    wire/.style={->, thick},
]

  \node[block] (G) at (0,0)   {$\Gamma$\\\scriptsize input-side\\[-2pt]\scriptsize coupling};
  \node[block] (D) at (2.2,0) {$\Delta$\\\scriptsize discarded\\[-2pt]\scriptsize dynamics};
  \node[block] (T) at (4.4,0) {$\Theta$\\\scriptsize output-side\\[-2pt]\scriptsize coupling};

  \draw[wire] (-2.15,0) -- node[above]{\scriptsize$u$} (G.west);
  \draw[wire] (G.east) -- (D.west);
  \draw[wire] (D.east) -- (T.west);
  \draw[wire] (T.east) -- node[above]{\scriptsize$y_e$} (6.55,0);

  \draw[dashed, rounded corners=4pt]
    (-1.25,-1.0) rectangle (5.65,1.0);

\end{tikzpicture}
\caption{Product (cascade) representation of the PROM error system
$\Sigma_e = \Theta(\Sigma,P)\,\Delta(\Sigma,P)\,\Gamma(\Sigma,P)$,
separating input-side coupling, discarded dynamics, and output-side coupling.}
\label{fig:error-system-forms}
\end{figure}

We observe that each of the three systems of the PROM product form
\eqref{eq:PROM error realization} implicitly captures the contribution
to the reduction error of the parts of $\tilde{\Sigma}$ in  \eqref{eq:full order system PQ realization}
left out in $\Sigma_{r}(\Sigma,P)$: the map $\Theta(\Sigma,P)$ captures $A_{PQ}$ and $CQ$, the map $\Gamma(\Sigma,P)$ captures $A_{QP}$
and $Q^{T}B$ and the map $\Delta(\Sigma,P)$ captures the full-order dynamics viewed through the discarded subspace.

\begin{cor}[Frequency-uniform product bound]\label{cor:productbound}
If $A$ and $A_{PP}$ are Hurwitz, then all three factors in \eqref{prod_form_elements} are stable and
\begin{equation}
\|\Sigma_e\|_\infty \leq
\underbrace{\|\Theta(\Sigma,P)\|_\infty\|\Delta(\Sigma,P)\|_\infty\|\Gamma(\Sigma,P)\|_\infty}_{:=b_{\rm prod}(\Sigma,P)}.
\label{eq:productbound}
\end{equation}
\end{cor}
\begin{proof}
The state matrices of the three factors are $A_{PP}$, $A$, and $A_{PP}$.  The result follows from submultiplicativity of the induced $\mathcal L_2$ norm.
\end{proof}

The certificate in \eqref{eq:productbound} is uniform over frequency.  In fact, for every $\omega\in\R$,
\begin{align}
\bar\sigma\!\left(\hat\Sigma_e(j\omega)\right)
&\leq
\bar\sigma\!\left(\hat\Theta(j\omega)\right)
\bar\sigma\!\left(\hat\Delta(j\omega)\right)
\bar\sigma\!\left(\hat\Gamma(j\omega)\right)\nonumber\\
&\leq b_{\rm prod}(\Sigma,P).\nonumber
\end{align}
Equivalently, under zero initial conditions, every $u\in\mathcal L_2$
satisfies the input--output uncertainty bound
\begin{equation*}
\|y-y_r\|_{\mathcal L_2}\leq
b_{\rm prod}(\Sigma,P)\,\|u\|_{\mathcal L_2}.
\end{equation*}
Thus a single scalar controls the error system for all frequencies and all
finite-energy inputs.  The condition on $A_{PP}$ is essential, since a Hurwitz
nonsymmetric $A$ can have an unstable orthogonal compression.  Moreover,
this is a certificate for a chosen projection, not an a priori bound that is
uniform over all admissible IIPROMs or depends only on the reduced order.
Without a common stability margin, admissible reduced poles can approach the
imaginary axis and the factor norms can become arbitrarily large.  A uniform
bound over a family of projections therefore requires common stability and
factor-norm bounds.  Finally, evaluating $b_{\rm prod}$ requires two
reduced-order norm computations and one full-order norm computation, so it
is not presented as a universal complexity reduction.  It instead decomposes
the error by mechanism and supplies a useful score when the factor norms or
their dc gains are inexpensive to evaluate.

A first observation related to the product form of Theorem \ref{thm:PROMbound} is that if the PROM is an IIPROM, then the realizations of $\Theta(\Sigma,P)$
and $\Gamma(\Sigma,P)$ are strictly-proper.

\begin{cor}[IIPROM error system]\label{cor_IIPROMproduct}
If  $\Sigma_{r}(P)$
is an IIPROM, then $\Theta(\Sigma,P)$ in \eqref{lti.1}
and $\Gamma(\Sigma,P)$ in \eqref{lti.2}
are strictly proper with realizations
\begin{equation}
\Theta(\Sigma,P):=\left(A_{PP},A_{PQ},CP,0_{n_y\times(n_x-r)}\right),
\end{equation}
and
\begin{equation}
\Gamma(\Sigma,P):=\left(A_{PP},P^T B,A_{QP},0_{(n_x-r)\times n_u}\right).
\end{equation}
\end{cor}
\begin{proof}
Since $\Sigma_{r}(P)$ is an
IIPROM, from Definition \ref{IIPROM} we have
$C=CPP^T$ and $B=PP^T B$. From $Q^T P=0$
and $Q^T Q=I_{n_x-r}$ we get $PP^T=I_{n_x}-Q Q^T$
such that $C=C\left(I_{n_x}-QQ^T\right)$ and
$B=\left(I_{n_x}-QQ^T\right)B$, therefore, $CQ=0$
 and $Q^T B=0$, obtaining the stated zero direct terms.
\end{proof}

Interface invariance removes an instantaneous error path through the discarded coordinates, so the error transfer matrix is strictly proper---there is no instantaneous feedthrough from input to error output.

\section{Exact Error for Positively Dominated Factors}\label{sec.positive}
The product bound~\eqref{eq:PROM error realization} is tight for an important special class of systems---the \emph{positively dominated} systems introduced by Rantzer~\cite{rantzer2012optimizing}.  We next recall the relevant definitions before deriving the main exactness result.

Positive LTI systems have the unique property of maintaining output within an invariant subspace $\R_+^p$ for non-negative input signals. They are formally defined below.

\begin{definition}[\cite{farina2000positive}]
An LTI system $\Sigma:=\left(A,B,C,D\right)$ is \emph{internal positive} if for every non-negative initial condition $x(0)$ and every non-negative input $u(t)$, $t \geq 0$,
the associated state and output evolutions $x(t)$ and
$y(t)$, $t \geq 0$, remain non-negative at every time instant
$t$.\end{definition}

Positive systems admit a special realization.
\begin{prop}[\cite{farina2000positive}]\label{Positive system realization} A system $\Sigma$ with realization $(A,B,C,D)$
is \emph{internally positive} if and only if $A$ is Metzler and $B,C,D\geq0$.
\end{prop}

Rantzer introduced \emph{positively dominated} systems as a generalization
of positive systems \cite{rantzer2012optimizing}.
\begin{definition}[Positively dominated system]
A system $\hat \Sigma(s)\in \mathcal{RH}_{\infty}^{n\times m}$
is called \emph{positively dominated} if every matrix entry satisfies
$\left|\hat{\Sigma}_{kl}\left(j\omega\right)\right|\leq\hat{\Sigma}_{kl}\left(0\right)$
for $\omega\in\R$. The set of all such matrices is denoted $P\mathcal{H}_{\infty}^{n\times m}$.
\end{definition}
Positively dominated systems have their $\mathcal H_{\infty}$ norm at dc and are closed under compatible products.  For internally positive systems the first property is the familiar fact that the $\mathcal H_\infty$ norm equals the dc gain \cite{tanaka2011bounded}.

\begin{prop}[\cite{rantzer2012optimizing}]\label{positively-dominated-properties}

    Let $\hat{\Sigma}\in P\mathcal{H}_{\infty}^{n\times m}$, then $\|\Sigma\|_{\mathcal H_{\infty}}=\|\hat{\Sigma}\left(0\right)\|_{2}$.
  Furthermore, for $\hat{\Sigma}_{1}\in P\mathcal{H}_{\infty}^{n\times m}$ and $\hat{\Sigma}_{2}\in P\mathcal{H}_{\infty}^{l\times n}$,
then $\hat{\Sigma}_{2}\hat{\Sigma}_{1}\in P\mathcal{H}_{\infty}^{l\times m}$.

\end{prop}

The tightness condition of Proposition~\ref{prop:sequence-tightness} takes a particularly simple form for positively dominated step errors: the peak frequency is pinned to dc, reducing the alignment condition to an algebraic criterion on the dc gain matrices.

\begin{cor}[Positive collocated steps]\label{cor:positive-sequence-tightness}
Suppose every step error $E_k$ is positively dominated and $E_k(0)$ is
symmetric positive semidefinite.  Then the sequence bound
\eqref{eq:Singleton PROM sequence bound} is tight if and only if the
dominant eigenspaces of all nonzero $E_k(0)$ share a common unit vector.
In particular, if
\begin{equation}
E_k(0)=\frac{z_kz_k^T}{d_k},\qquad d_k>0,
\label{eq:rank-one-step-error}
\end{equation}
then equality holds if and only if all nonzero $z_k$ are collinear.  The
condition is automatic for scalar input and output.
\end{cor}
\begin{proof}
Positive dominance pins the common peak frequency to $\omega_\star=0$.  For
a positive semidefinite dc matrix, the spectral norm is its largest
eigenvalue, so \eqref{eq:sequence-alignment} is equivalent to a common
dominant eigenvector.  The dominant eigenspace of the rank-one matrix in
\eqref{eq:rank-one-step-error} is $\operatorname{span}(z_k)$.
\end{proof}

\begin{remark}[Order invariance for rank-one positive steps]
\label{rem:sequence-order-invariance}
Under the rank-one hypotheses of Corollary~\ref{cor:positive-sequence-tightness},
positive dominance gives $\|E_k\|_\infty=\|E_k(0)\|_2=\operatorname{trace}E_k(0)$
for each step.  Summing over $k$, the right-hand side of \eqref{eq:Singleton PROM sequence bound}, $U_r:=\sum_k\|E_k\|_\infty$, equals $\operatorname{trace}E(0)$, where
$E(0)=\sum_kE_k(0)$ is the final dc error---independent of deletion order.
Consequently, reordering cannot improve the bound, and it is tight if and
only if $\operatorname{rank}E(0)=1$.
\end{remark}

\begin{thm}[Exact error at dc]\label{thm:positive-error}
Let $A$ and $A_{PP}$ be Hurwitz.  If the three factors $\Theta(\Sigma,P)$, $\Delta(\Sigma,P)$, and $\Gamma(\Sigma,P)$ are positively dominated, then $\Sigma_e$ is positively dominated and
\begin{equation}
\|\Sigma_e\|_\infty=\|\hat\Sigma_e(0)\|_2
=\left\|CA^{-1}Q(Q^TA^{-1}Q)^{-1}Q^TA^{-1}B\right\|_2.
\label{eq:positive-exact}
\end{equation}
\end{thm}
\begin{proof}
Closure under multiplication and Theorem \ref{thm:PROMbound} imply
that $\Sigma_e$ is positively dominated.  Proposition
\ref{positively-dominated-properties} therefore gives $\|\Sigma_e\|_\infty=\|\hat\Sigma_e(0)\|_2$. Substitution of $s=0$ in
\eqref{eq:PROM error TFM} gives \eqref{eq:positive-exact}; the four minus
signs cancel.
\end{proof}

The hypotheses concern the factors, not merely the original realization.  Positivity of $\Sigma$ alone does not guarantee positivity of the factors for an arbitrary dense orthogonal basis because $P$ and $Q$ may have mixed signs.  A simple sufficient case is coordinate selection: if $\Sigma$ is stable and internally positive and $P=e_{\mathcal I}^n$, then $Q=e_{\mathcal I^c}^n$ (up to permutation), $A_{PP}$ is a Hurwitz principal submatrix of the Hurwitz Metzler matrix $A$, and all factor realizations are internally positive.  Equation \eqref{eq:positive-exact} then applies.

The exact dc formula also yields an a priori upper bound under a
strict dissipativity assumption.  The setting is deliberately narrower than
Theorem~\ref{thm:PROMbound}: the projection must be interface invariant,
$-A$ must be uniformly coercive, and positive dominance is needed to
promote the dc estimate to a uniform $\mathcal{H}_\infty$ certificate.

To state the result compactly, let $H:=-A$, complete $P$ to an orthogonal
matrix $[P\ Q]$, and set
\begin{align}
H_{PP}&:=P^THP, & H_{PQ}&:=P^THQ,\nonumber\\
H_{QP}&:=Q^THP, & H_{QQ}&:=Q^THQ.\nonumber
\end{align}
Define the Schur complement $S_Q:=H_{QQ}-H_{QP}H_{PP}^{-1}H_{PQ}$ and the
coupling factors
\begin{equation*}
L_P:=CPH_{PP}^{-1}H_{PQ},\qquad R_P:=H_{QP}H_{PP}^{-1}P^TB.
\end{equation*}

\begin{thm}[Coercive error certificate]\label{thm:coercive-error}
Let $\Sigma=(A,B,C,D)$ and let $P$ define an IIPROM.  Suppose
\begin{equation}
\tfrac{1}{2}(H+H^T)\succeq \eta I_{n_x}
\label{eq:H-coercivity}
\end{equation}
for some $\eta>0$.  Then $H_{PP}$ and $S_Q$ are nonsingular, and
\begin{align}
\hat\Sigma_e(0)&=L_P S_Q^{-1}R_P,
\label{eq:dc-schur-error}\\
\|\hat\Sigma_e(0)\|_2
&\leq \frac{\|L_P\|_2\|R_P\|_2}{\eta}
\leq
\frac{\|C\|_2\|B\|_2\|H_{PQ}\|_2\|H_{QP}\|_2}{\eta^3}.
\label{eq:coercive-dc-bound}
\end{align}
If in addition $\Sigma_e$ is positively dominated---in particular if the three
factors in Theorem~\ref{thm:positive-error} are---then
\begin{equation}
\|\Sigma_e\|_\infty=\|\hat\Sigma_e(0)\|_2\leq \frac{\|L_P\|_2\|R_P\|_2}{\eta}.
\label{eq:coercive-Hinf-bound}
\end{equation}
\end{thm}
\begin{proof}
Interface invariance gives $Q^TB=0$ and $CQ=0$.  Applying the block-inverse
formula to $[P\ Q]^TH[P\ Q]$ and subtracting the reduced dc map
$CPH_{PP}^{-1}P^TB$ from the full dc map $CH^{-1}B$ gives
\eqref{eq:dc-schur-error}.  Condition~\eqref{eq:H-coercivity} implies
$(H_{PP}+H_{PP}^T)/2\succeq\eta I_r$, and hence
$\|H_{PP}^{-1}\|_2\leq1/\eta$.  For any $y$, choose
$x=-H_{PP}^{-1}H_{PQ}y$.  Then
\begin{equation*}
\operatorname{Re}\!\left(y^*S_Qy\right)
=\operatorname{Re}\!\left(
\begin{bmatrix}x\\y\end{bmatrix}^{\!*}
[P\ Q]^TH[P\ Q]
\begin{bmatrix}x\\y\end{bmatrix}
\right)
\geq\eta\|y\|_2^2.
\end{equation*}
Thus $\|S_Q^{-1}\|_2\leq1/\eta$.  Submultiplicativity gives the first
inequality in \eqref{eq:coercive-dc-bound}; the second follows from
$\|P\|_2=\|Q\|_2=1$ and the bound on $H_{PP}^{-1}$.  Positive dominance
gives \eqref{eq:coercive-Hinf-bound} by Proposition
\ref{positively-dominated-properties}.
\end{proof}

\begin{remark}[Collocated symmetric case]\label{rem:collocated-coercive}
If additionally $H=H^T$ and $C=B^T$, set $Z_P:=H_{QP}H_{PP}^{-1}P^TB=L_P^T$.
Then $S_Q\succ0$, $\hat\Sigma_e(0)=Z_P^TS_Q^{-1}Z_P\succeq0$, and under
positive dominance $\|\Sigma_e\|_\infty\leq\|Z_P\|_2^2/\eta$.  Indeed,
symmetry gives $H_{PQ}=H_{QP}^T$ and $H_{PP}=H_{PP}^T$, so
$L_P=CPH_{PP}^{-1}H_{PQ}=(H_{QP}H_{PP}^{-1}P^TB)^T=Z_P^T$ and
$R_P=Z_P$; $S_Q$ is a Schur complement of the positive-definite matrix
$[P\ Q]^TH[P\ Q]$ and is therefore positive definite.
\end{remark}

The word ``a priori'' here is projection dependent. Once a candidate
$P$ is specified, \eqref{eq:coercive-Hinf-bound} can be evaluated without
forming the full error system, computing $S_Q^{-1}$, or performing a frequency
sweep.  Partition-dependent a priori bounds of a similar kind
are derived for clustering of leader--follower multi-agent systems in
\cite{jongsma2018model}.  The bound \eqref{eq:coercive-Hinf-bound} does not depend only on the target order, and it is not an
approximation ratio relative to the globally optimal projection.  Coercivity
alone bounds only the dc error. Without positive dominance, an off-zero
frequency peak can exceed \eqref{eq:coercive-dc-bound}.  The coordinate-selection
specialization of a stable internally positive system satisfies the needed
positive-dominance condition, whereas a general dense projection or graph
aggregation need not.

\section{Model Reduction of MAS\\ Using Graph Contractions}\label{sec.MAS_IIPROM}

The preceding sections characterize admissible projections and their error, but do not prescribe how to choose the free subspace in Theorem~\ref{thm:IIPROM-geometry}.  For multi-agent systems it is natural to require every reduced coordinate to represent one original agent or a connected aggregate, with actuated and measured agents retaining individual identities.  Vertex partitions provide exactly this structure: their normalized characteristic matrices define orthogonal projections whose quotient graph gives a physical interpretation to the reduced dynamics.  Requiring every interface vertex to be a singleton cell enforces the protected-subspace condition of Theorem~\ref{thm:IIPROM-geometry}, replacing the continuous search over projection subspaces with a combinatorial search over interpretable quotient graphs.  When the positively dominated specialization of Section~\ref{sec.positive} applies, the true $\mathcal{H}_\infty$ error equals the dc gain, so candidate contractions can be scored without a frequency sweep.  This construction is a form of graph coarsening \cite{loukas2019graph} and is related to clustering reductions \cite{ishizaki2014model,monshizadeh2014projection,besselink2016clustering}, in particular to the positivity-preserving clustering of positive directed networks \cite{ishizaki2015clustered} and to dissimilarity-based clustering \cite{cheng2019dissimilarity}; its distinctive constraint is exact preservation of every interface vertex.

\subsection{Edge-Matching Graph Partitions}\label{subsec.edgecontract}

\begin{figure}[!h]
    \centering
    \subfloat[]{\includegraphics[width=0.3\columnwidth]{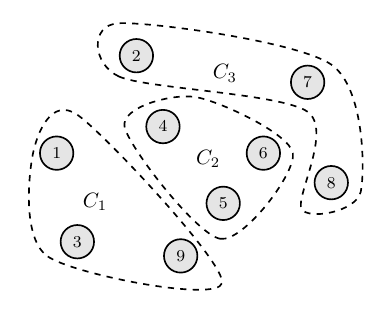}\label{fig:rpart}}\hfill
    \subfloat[]{\includegraphics[width=0.3\columnwidth]{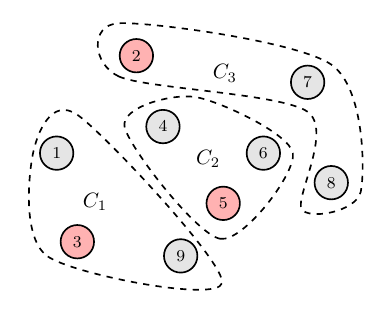}\label{fig:Iinv}}\hfill
    \subfloat[]{\includegraphics[width=0.3\columnwidth]{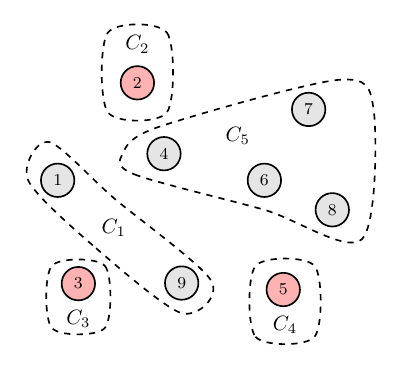}\label{fig:strongIinv}}
    \caption{A set of $|\V|=9$ nodes showing (a)~a 3-partition, (b)~an $\mathcal I$-invariant 3-partition with $\mathcal I=\{2,3,5\}$ (in red), and (c)~a strongly $\mathcal I$-invariant 5-partition.}
    \label{fig:partitions}
\end{figure}

An $r$-\emph{partition}, denoted as $\pi_r$, of a vertex set $\V$ is the partition $\{C_i\}_{i=1}^r$ of $\V$ into $r$ cells such that
$$\bigcup_{i=1}^r C_i = \V \text{ and } C_i \cap C_j = \emptyset.$$
We call an $r$-partition $\mathcal I$-\emph{invariant} for some $\mathcal I \subset \V$ if for each $i=1,\ldots,r$, one has $|C_i \cap \mathcal I| \leq 1$.  It is \emph{strongly} $\mathcal I$-invariant if
\begin{equation}
C_i\cap\mathcal I\neq\emptyset\quad\Longrightarrow\quad |C_i|=1.
\label{eq:strongpartition}
\end{equation}
Thus every interface vertex forms a singleton cell. We denote the set of strongly $\mathcal I$-invariant $r$-partitions by $S_r(\V)$.  Figure \ref{fig:partitions} illustrates the definitions.

\begin{figure}[!h]
    \centering
    \includegraphics[width=0.18\textwidth]{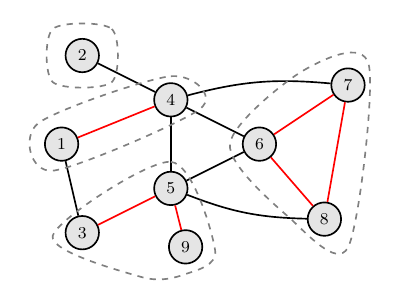}
    \caption{An edge-induced partition of a connected graph $\G$. The edges in $\E_s$ are drawn in red and the resulting connected components of the induced graph $\hat \G = (\V,\E_s)$ are outlined with dashed lines. Thus the edge-induced partition $\pi(\E_s)$ yields a $4$-partition.}
    \label{fig:edgepart}
\end{figure}

Given a graph $\G = (\V, \E)$ and a subset of its edges $\E_s \subseteq \E$, the \emph{edge-induced partition} $\pi(\E_s)$ is the collection of connected components of $(\V,\E_s)$ \cite{leiter2021edge}. This is illustrated in Figure \ref{fig:edgepart}.

Given an $r$-partition $\pi = \{C_i\}_{i=1}^r$, we define $M(\pi) \in \R^{n \times r}$ as the \emph{partition characteristic matrix} (PCM). Its entry $[M(\pi)]_{ij}$ equals $1$ when $i \in C_j$ and $0$ otherwise.  The corresponding
partition projection matrix (PPM) is $P(\pi)=M(\pi)(M(\pi)^TM(\pi))^{-1/2}$. For example, for the $3$-partition in Figure \ref{fig:rpart}, we have
$$M(\pi) = \begin{bmatrix} 1 & 0 & 1 & 0 & 0 & 0 & 0 & 0 &  1\\0 & 0 & 0 & 1 & 1 & 1 & 0 & 0 & 0 \\0 & 1 & 0 & 0 & 0 & 0 & 1 & 1 & 0\end{bmatrix}^T.$$

For a MAS $\Sigma(\G,\U,\Y)$ let $\mathcal I=\U\cup\Y$.  If $B$ and $C$ select interface coordinates and $\pi\in S_r(\V)$, then $P(\pi)P(\pi)^T$ fixes each interface basis vector.  Consequently $B=P(\pi)P(\pi)^TB$ and $C=CP(\pi)P(\pi)^T$, so the partition PROM is interface invariant.  Strong rather than ordinary invariance is required. Merging an interface vertex with a non-interface vertex changes its coordinate vector.

\subsection{Stability of Positive MAS Contractions}\label{subsec:positive-contraction-stability}

For Metzler network dynamics, graph contractions admit a stability certificate that is distinct from the Euclidean dissipativity condition of Proposition~\ref{prop:compression-stability}.  Recall that a Metzler matrix $A$ is Hurwitz if and only if there exists a vector $v>0$ such that $Av<0$.  A partition projection preserves this certificate whenever the positive vector is constant on every partition cell.

\begin{prop}[Positive-vector stability of MAS contractions]\label{prop:positive-vector-contraction}
Let $A\in\R^{n\times n}$ be Metzler, and suppose there exists $v>0$ such that $Av<0$.  Let $\pi=\{C_i\}_{i=1}^r$ be a partition with normalized characteristic projection $P=P(\pi)=M(\pi)(M(\pi)^TM(\pi))^{-1/2}$.  If $v\in\operatorname{im}P$, equivalently, $v$ is constant on every cell $C_i$, then the contracted matrix $A_r=P^TAP$ is Metzler and Hurwitz.  More precisely, writing $v=Pv_r$, one has $v_r>0$ and $A_rv_r=P^TAv<0$.
\end{prop}
\begin{proof}
Because $P$ is a normalized characteristic matrix, $P\geq0$ and its columns have disjoint supports.  The assumption $v\in\operatorname{im}P$ and $v>0$ therefore imply that the unique vector $v_r=P^Tv$ satisfying $v=Pv_r$ is strictly positive.  For $i\neq j$, $[A_r]_{ij}=p_i^TAp_j\geq0$, since $p_i,p_j\geq0$ have disjoint supports and all entries of $A$ connecting distinct coordinates are nonnegative; hence $A_r$ is Metzler.  Moreover, $A_rv_r=P^TAPv_r=P^TAv<0$, because $Av<0$ and every row of $P^T$ is nonnegative and nonzero.  The Metzler stability criterion then implies that $A_r$ is Hurwitz.
\end{proof}

\begin{remark}[Grounded-Laplacian networks]\label{rem:grounded-laplacian-contraction}
For $A=-(L+gI)$ with $g>0$ and a graph Laplacian satisfying $L\mathbf 1=0$, choose $v=\mathbf 1$.  Then $Av=-g\mathbf 1<0$, and $\mathbf 1$ is constant on every partition cell, so Proposition~\ref{prop:positive-vector-contraction} applies to every partition.  If $n_i=|C_i|$ and $v_r=[\sqrt{n_1},\ldots,\sqrt{n_r}]^T$, then $Pv_r=\mathbf 1$ and $A_rv_r=-g v_r$.
Thus every graph contraction is Hurwitz. Since $A_r$ is Metzler and has the strictly positive eigenvector $v_r$, its spectral abscissa is $-g$.  Hence the contraction preserves the decay margin exactly.
\end{remark}

This criterion is complementary to Proposition~\ref{prop:compression-stability}: it is specific to nonnegative partition projections and Metzler dynamics, but does not require Euclidean dissipativity and therefore also applies to nonsymmetric positive network models.  Positive dominance itself remains an input--output property of an already stable transfer matrix; Proposition~\ref{prop:positive-vector-contraction} supplies the internal stability certificate needed for this class of MAS contractions.

The natural goal is to find the partition $\pi^*\in S_r(\mathcal{V})$ that minimizes $\|\Sigma-\Sigma_r(P(\pi))\|_\infty$. Finding such an optimal partition, however, is numerically intractable for a moderate number of nodes: the number of $r$-partitions is the Stirling number of the second kind,
which for $r\ll n$ is asymptotically $S\left(n,r\right)\sim\frac{r^{n}}{r!}$
\cite[p.18]{wilf1994generatingfunctionology}.

\subsection{Sequential Greedy Coarsening}

Since exact minimization over $S_r(\mathcal{V})$ is combinatorially intractable, we adopt a greedy heuristic that builds the target partition one vertex-merge at a time. A naive approach would select the $n-r$ best edges simultaneously, but this can create cycles and need not produce exactly $r$ cells; it also ignores interactions between successive merges. Algorithm~\ref{greedy-edge-IIPROM} below addresses both issues: it performs one admissible merge at a time, recomputing every candidate score after each step. The default score is the true error $J(P)=\|\Sigma-\Sigma_r(P)\|_\infty$, directly addressing cases in which the product bound is conservative.  When exact norm evaluation is too costly, $b_{\rm prod}$ in \eqref{eq:productbound} may be substituted, but only after the stability tests required by Corollary~\ref{cor:productbound}.

The procedure is formalized in Algorithm~\ref{greedy-edge-IIPROM}. Starting from the finest partition (each vertex its own cell), at each iteration all edge-adjacent cell pairs are enumerated; merges that would violate the strong interface-invariance constraint~\eqref{eq:strongpartition} or yield an unstable reduced model are discarded, and the remaining candidates are scored.  The algorithm commits to the merge with the best score and repeats until the target order $r$ is reached or no admissible merge remains.

\begin{center}
\begin{algorithm}[h]
\caption{Sequential interface-preserving edge coarsening} \label{greedy-edge-IIPROM}

\textbf{Input:} MAS $\Sigma(\G,\U,\Y)$ of order $n$, candidate edges $\E_c\subseteq\E$, target order $r\geq|\mathcal{I}|$, score $J(\pi'):=\|\Sigma-\Sigma_{|\pi'|}(P(\pi'))\|_\infty$ or $b_{\rm prod}(\Sigma,P(\pi'))$, optional Lyapunov $X$.
\begin{enumerate}
\item \textbf{Initialize.} Set $\pi\leftarrow\{\{v\}:v\in\V\}$.  Test \eqref{eq:euclidean-dissipativity}; record $\mu$ if satisfied, else validate $X$.  For Metzler dynamics, check whether $v>0$ with $Av<0$ is available (Prop.~\ref{prop:positive-vector-contraction}).
\item \textbf{Repeat while} $|\pi|>r$\textbf{:}
  \begin{enumerate}
    \item[(a)] \textit{Generate.} Form all $\pi'$ by merging adjacent cells via edges in $\E_c$.
    \item[(b)] \textit{Filter.} Discard $\pi'$ violating \eqref{eq:strongpartition}.  Certify stability via Prop.~\ref{prop:compression-stability}; for Metzler dynamics use Prop.~\ref{prop:positive-vector-contraction} when $v\in\operatorname{im}P(\pi')$.  Remaining uncertified candidates undergo a direct Hurwitz test; unstable ones are discarded.  If $J=b_{\rm prod}$, also discard candidates with unstable factors.
    \item[(c)] \textit{Select.} If no admissible $\pi'$ remains, return \texttt{early-stop}.  Otherwise set $\pi\leftarrow\operatorname{argmin}_{\pi'} J(\pi')$ and record $J(\pi)$, $b_{\rm prod}$, and the certificate.
  \end{enumerate}
\end{enumerate}
\textbf{Output:} $\Sigma_{|\pi|}(P(\pi))$, termination status, and per-step scores and certificates.
\end{algorithm}
\par\end{center}

The stability certificates used in Algorithm~\ref{greedy-edge-IIPROM} are Proposition~\ref{prop:compression-stability} and, for positive network contractions, the certificate of Proposition~\ref{prop:positive-vector-contraction} tailored to nonnegative partition projections.

Taking $\E_c=\E$ gives the richest local search and is the default in the supplied implementation.  The merge objective is not known to be submodular, so data-independent greedy approximation ratios do not apply. The following result instead provides a computable a posteriori certificate for every run that reaches the requested order.

\begin{prop}[A posteriori greedy certificate]\label{prop:greedy-certificate}
Let $\sigma_{r+1}^{\rm H}(\Sigma)$ be the $(r+1)$st Hankel singular value of the stable full-order system, with value zero when its McMillan degree is at most $r$.  Suppose Algorithm~\ref{greedy-edge-IIPROM} reaches the requested order with a stable sequence $\{\Sigma_{(k)}^{\rm g}\}_{k=0}^{n-r}$, and define
\begin{equation*}
J_r^{\rm g}:=\|\Sigma-\Sigma_{(n-r)}^{\rm g}\|_\infty,
\qquad
U_r^{\rm g}:=\sum_{k=1}^{n-r}
\|\Sigma_{(k)}^{\rm g}-\Sigma_{(k-1)}^{\rm g}\|_\infty.
\end{equation*}
Let $J_r^\star$ be the optimal value of Problem~\ref{prob:optimalIIPROM}, and let $V_r^\star$ be the optimal value of Problem~\ref{prob.sequencebound}.  Then
\begin{equation}
\sigma_{r+1}^{\rm H}(\Sigma)\leq J_r^\star\leq J_r^{\rm g}\leq U_r^{\rm g},
\qquad
\sigma_{r+1}^{\rm H}(\Sigma)\leq V_r^\star\leq U_r^{\rm g}.
\label{eq:greedy-certificate-chain}
\end{equation}
Consequently, if $\sigma_{r+1}^{\rm H}(\Sigma)>0$, the returned model satisfies the certified approximation factors
\begin{equation}
\frac{J_r^{\rm g}}{J_r^\star}
\leq\frac{J_r^{\rm g}}{\sigma_{r+1}^{\rm H}(\Sigma)},
\qquad
\frac{U_r^{\rm g}}{V_r^\star}
\leq\frac{U_r^{\rm g}}{\sigma_{r+1}^{\rm H}(\Sigma)}.
\label{eq:greedy-factor-certificate}
\end{equation}
\end{prop}
\begin{proof}
For every stable transfer matrix $\widetilde\Sigma_r$ of McMillan degree at most $r$, the Hankel-operator approximation theorem gives
\begin{equation*}
\|\Sigma-\widetilde\Sigma_r\|_\infty
\geq \|\mathcal H_\Sigma-\mathcal H_{\widetilde\Sigma_r}\|
\geq \sigma_{r+1}^{\rm H}(\Sigma)
\end{equation*}
\cite{glover1984all,antoulas2005approximation}.  This applies to every feasible final model in both problems.  The triangle inequality gives $J_r^{\rm g}\leq U_r^{\rm g}$ and, for every feasible sequence, its accumulated cost is no smaller than its final error.  Finally, the greedy output is feasible for both problems whenever it terminates.  These observations give \eqref{eq:greedy-certificate-chain}; division by a positive lower bound gives \eqref{eq:greedy-factor-certificate}.
\end{proof}

Proposition~\ref{prop:greedy-certificate} is an a posteriori result; it does not establish a uniform approximation ratio.  Algorithm~\ref{greedy-edge-IIPROM} returns the certified upper bound $U_r^{\rm g}$, while the Hankel singular value supplies a universal lower bound on either global optimum.  The lower bound is taken over all stable order-$r$ approximants, including dense unstructured ones, and can be much smaller than the best interface-preserving error when the plant is nearly reducible in unstructured coordinates---as in the thermal network of Section~\ref{sec.cases}, where the interface-preserving error exceeds $\sigma^{\rm H}_{r+1}$ by a factor of at least $277$ already at $r=60$. The certificate is therefore a worst-case safety net rather than a typical performance prediction.

\section{Case Study: A $65{,}566$-State Chip Thermal Network}\label{sec.cases}

Dynamic thermal management of a processor acts on the power of individual functional units and reads their temperatures, so a compact thermal model is useful only if every unit keeps its own input and output.  The network below is internally positive, symmetric, and strictly dissipative, so every result of Sections~\ref{sec.problem}--\ref{sec.MAS_IIPROM} applies to it.  The structured baseline is Kron reduction of the same interface, known in thermal modeling as static condensation or Guyan reduction \cite{guyan1965reduction}: it keeps the interface nodes and $r-|\mathcal I|$ interior nodes and eliminates the rest through a Schur complement.  It is a Petrov--Galerkin PROM that preserves the interface and matches the dc gain exactly, but it is not an orthogonal projection, so Sections~\ref{sec.structure}--\ref{sec.positive} do not apply to it.  We also compare against coordinate deletion and balanced truncation (BT), and report the true relative error $\|\Sigma-\Sigma_r\|_\infty/\|\Sigma\|_\infty$.

\subsection{Model and Reduction Methods}

The network follows the HotSpot methodology \cite{huang2006hotspot}, but not its solver, for the public Alpha EV6 floorplan and GCC power trace.  A $150\,\mu$m silicon die ($k=130$\,W/(m\,K), $1.63$\,MJ/(m$^3$\,K)) is discretized into a $256\times256$ grid.  Each of the $30$ functional units has one active-layer node ($10\,\mu$m) coupled to the cells beneath it through a $20\,\mu$m interface layer ($4$\,W/(m\,K)), and the grid is cooled through a uniform $0.1$\,K/W package resistance, giving $n=65{,}566$.  With conductance matrix $G$ (a grounded Laplacian), diagonal heat capacities $C_{\rm th}$, unit-node selector $S$, and energy coordinates $x=C_{\rm th}^{1/2}(T-T_{\rm amb})$, the model is
\begin{equation}
\dot x=-Hx+Bu,\qquad y=B^Tx,
\label{eq:hotspot-model}
\end{equation}
with $H=C_{\rm th}^{-1/2}GC_{\rm th}^{-1/2}$ and $B=C_{\rm th}^{-1/2}S$; the inputs are the unit powers and the outputs are the unit temperature rises.  Here $H\succeq\eta I$ with $\eta=149.6$\,s$^{-1}$, $\|\Sigma\|_\infty=\|S^TG^{-1}S\|_2=44.9$\,K/W, and $r_{\min}=30$ by Theorem~\ref{thm:IIPROM-geometry}.

\emph{IIPROM.}  Every unit node is a singleton cell, and $P(\pi)=C_{\rm th}^{1/2}M(M^TC_{\rm th}M)^{-1/2}$, $M=M(\pi)$, is the capacitance-weighted partition projection of Section~\ref{sec.MAS_IIPROM}: it is orthonormal and nonnegative and contains $\mathcal S_{\rm io}$.  In temperature coordinates the reduced model is an RC network on the quotient graph whose states are cluster-mean temperatures.  The partition is built by a one-pass divisive heuristic rather than by Algorithm~\ref{greedy-edge-IIPROM}.  Starting from one cluster per unit, it repeatedly splits into quadrants the cluster with the largest capacitance-weighted spread $\sum_{i\in C}c_i\|T_i-\bar T_C\|_2^2$ of the dc temperature field $T=G^{-1}S$; these spreads sum to $\|(I-PP^T)X_0\|_F^2$, $X_0=H^{-1}B$, the dc response discarded by aggregation.  Cells never cross unit boundaries, so the smallest order is $r=60$.

\emph{Baselines.}  Kron reduction retains either the $r-30$ cells of highest dc-influence $\|T_i\|_2$ (Kron-dc) or one representative cell per IIPROM cluster (Kron-rep).  Coordinate deletion keeps the Kron-rep nodes but truncates orthogonally, $P=e_{\mathcal I}$, the setting of Theorem~\ref{thm:positive-error}.  Since $H$ is symmetric and $C=B^T$, BT is the orthogonal PROM onto the $r$ leading Gramian eigenvectors, whose range does not contain $\mathcal S_{\rm io}$.  The Gramian factor is computed by low-rank ADI \cite{penzl2000cyclic} (rank $1440$, $20.6$\,s), so BT is reported for $r\le801$.  All reduced models are stable with decay margin at least $\eta$, the orthogonal ones by Proposition~\ref{prop:compression-stability}(i) and Kron because a Schur complement of $H$ satisfies the same bound.  Errors are computed from sparse factorizations of the full model at $90$ log-spaced frequencies and at $\omega=0$, refining every off-dc peak.

\begin{figure*}[t]
\centering
\includegraphics[width=\textwidth]{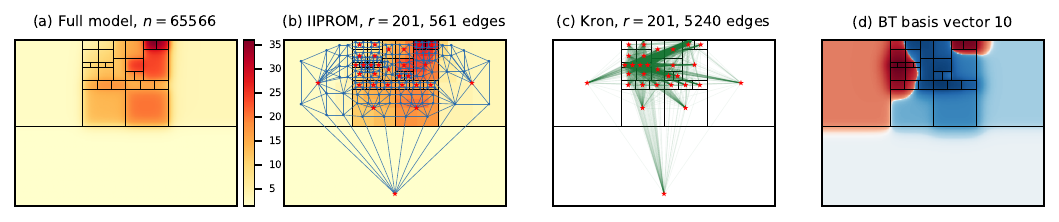}
\caption{EV6 thermal network and reduced models at $r=201$, with the uniform lower part of the L2 cache cropped; red stars mark the $30$ unit (interface) nodes.  (a)~Steady-state temperature rise (K) under the mean GCC power ($40.2$\,W).  (b)~IIPROM: reduced graph (blue) over the cluster-mean temperatures.  (c)~Kron-dc: retained cells and their $5{,}240$-edge reduced graph.  (d)~A BT basis vector (symmetric log scale), which corresponds to no physical location.}
\label{fig:hotspot-network}
\end{figure*}

\begin{table}[t]
\centering
\caption{EV6 thermal network: accuracy, structure, and cost}
\label{tab:hotspot}
\footnotesize
\setlength{\tabcolsep}{3pt}
\begin{tabular}{@{}llrrrccrr@{}}
\hline
$r$ & Method & Err.\,(\%) & $\omega_{\rm pk}$ & Edges & Pos. & Int. & Cap. & Time\\
\hline
402 & IIPROM   & $24.0$ & $0$   & $1{,}199$  & \checkmark & \checkmark & $100\%$ & $0.06$\\
    & Kron-dc  & $34.0$ & $466$ & $12{,}144$ & \checkmark & \checkmark & $6.8\%$ & $1.2$\\
    & Kron-rep & $42.3$ & $633$ & $51{,}644$ & \checkmark & \checkmark & $6.8\%$ & $1.3$\\
    & Deletion & $58.5$ & $0$   & $372$      & \checkmark & \checkmark & $6.8\%$ & $<0.01$\\
    & BT       & $10^{-10}$ & $0.6$ & $80{,}601$ & $\times$ & $\times$ & -- & $20.6$\\
\hline
6402 & IIPROM   & $1.38$ & $0$   & $19{,}702$  & \checkmark & \checkmark & $100\%$ & $0.15$\\
     & Kron-dc  & $3.10$ & $406$ & $31{,}902$  & \checkmark & \checkmark & $15.4\%$ & $32.6$\\
     & Kron-rep & $5.13$ & $209$ & $284{,}305$ & \checkmark & \checkmark & $15.4\%$ & $22.7$\\
     & Deletion & $47.0$ & $0$   & $14{,}892$  & \checkmark & \checkmark & $15.4\%$ & $<0.01$\\
\hline
\end{tabular}

\vspace{2pt}
\parbox{\linewidth}{\scriptsize Pos.: $A_r$ Metzler and $B_r\geq0$.  Int.: interface exactly preserved.  Cap.: fraction of heat capacity retained.  Time: seconds, excluding the $0.46$\,s dc factorization shared by IIPROM and Kron-dc.}
\end{table}

\subsection{Results}

\emph{Structure.}  The IIPROM is a sparse lumped RC network (Fig.~\ref{fig:hotspot-network}, Table~\ref{tab:hotspot}).  It has about three edges per state at every order and retains all of the heat capacity, and each of its states is the mean temperature of an identifiable region of the die.  The Kron-reduced model is also an RC network on physical nodes, because a Schur complement of an M-matrix is an M-matrix.  However, it discards $85$--$94\%$ of the heat capacity and fills in: at $r=402$, Kron-dc has $10$ times and Kron-rep $43$ times as many edges as the IIPROM.  BT yields a complete graph with a non-Metzler state matrix and a mixed-sign input matrix.  Its range only approximates the interface ($\|(I-PP^T)B\|_2/\|B\|_2=4.7\%$ at $r=60$), and its coordinates are global modes (Fig.~\ref{fig:hotspot-network}d).

\begin{figure}[t]
\centering
\includegraphics[width=0.95\linewidth]{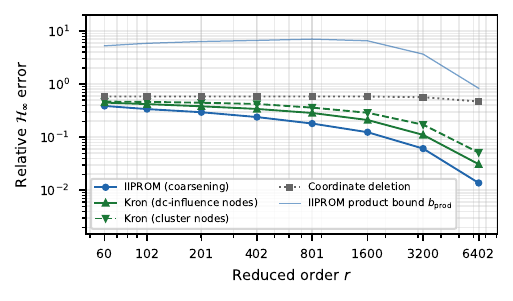}
\caption{EV6 thermal network: relative $\mathcal{H}_\infty$ error of the structure-preserving reductions vs.\ reduced order, with the IIPROM product bound of Corollary~\ref{cor:productbound}.}
\label{fig:hotspot-error}
\end{figure}

\emph{Accuracy.}  Among the structure-preserving reductions (Fig.~\ref{fig:hotspot-error}), the IIPROM is the most accurate at every order from $r=60$ to $r=6{,}402$.  The better Kron variant has $1.16$ times the IIPROM error at $r=60$ ($44.9\%$ vs.\ $38.8\%$), $1.6$ times at $r=801$, and $2.25$ times at $r=6{,}402$ ($3.10\%$ vs.\ $1.38\%$).  Coordinate deletion stays at $47$--$58\%$, because deleting a cell grounds it to ambient.  Per edge the advantage widens: the IIPROM reaches $18\%$ error with $2{,}466$ edges and $6.1\%$ with $9{,}981$, whereas Kron-dc needs $22{,}036$ and $35{,}913$ edges for $28.6\%$ and $11.1\%$.  BT is far more accurate, with $0.29\%$ at $r=60$ and machine precision from $r=402$, and it tracks the Hankel lower bound within a factor of about $2.5$.  The full model is therefore almost exactly reducible in unstructured coordinates, and the gap between BT and every node-based model is the price of keeping physical structure.  For the same reason, the a posteriori factor of Proposition~\ref{prop:greedy-certificate} is uninformative here: $J_r/\sigma_{r+1}^{\rm H}\geq277$ already at $r=60$.

\emph{Where the error lives.}  The IIPROM and deletion errors peak at $\omega=0$ at every tested order (Table~\ref{tab:hotspot}).  For deletion this is guaranteed by Theorem~\ref{thm:positive-error}.  For the IIPROM it is an observation rather than a certified property, since $Q$ has mixed signs, but it means that the dc-gain error, obtained from a single sparse solve, equals the $\mathcal H_\infty$ error.  Kron reduction is exact at dc, and its error peaks at $209$--$812$\,rad/s, above $\eta$, where the discarded heat capacity matters.  Kron therefore gives exact steady-state temperatures for power profiles slower than about $1/\eta=6.7$\,ms, whereas the IIPROM has the smaller worst-case error over all frequencies.

\emph{Certificates and cost.}  The product bound $b_{\rm prod}$ holds at every order, but it is conservative by a factor of $14$ at $r=60$ and $60$ at $r=6{,}402$, because it ignores how the input and output couplings align.  The coercive bound of Remark~\ref{rem:collocated-coercive} exceeds $\|\Sigma\|_\infty$ and is uninformative here.  Once the dc field is available ($0.46$\,s), an IIPROM is built in $0.02$--$0.15$\,s by aggregation alone.  Kron reduction takes $0.3$\,s at $r=60$ and $33$\,s at $r=6{,}402$, and BT needs the Gramian and yields dense models.  The complete study runs in $16$ minutes on a laptop.

\section{Conclusions}\label{sec.conclusion}

We characterized interface-invariant orthogonal projections as those whose
retained subspace contains $\operatorname{im}[B\ C^T]$, and factored the
reduction error into input-coupling, internal-dynamics, and output-coupling
systems.  The factorization gives a submultiplicative $\mathcal{H}_\infty$
bound that becomes an exact dc-gain expression for positively dominated
factors; strict dissipativity certifies stability, and a Hankel singular value
gives an a posteriori approximation factor.

On a $65{,}566$-state processor thermal network with $30$ functional-unit
interfaces, interface-preserving graph coarsening produced sparse, positive RC
models that retain all of the heat capacity.  At every tested order they were
more accurate in $\mathcal{H}_\infty$ than Kron reduction of the same
interface, by up to a factor of $2.25$ and with up to $43$ times fewer edges.
Balanced truncation was far more accurate but produced dense models whose
coordinates mix the interface; this gap is the price of keeping physical
structure.

\appendices
\section{Proof of Theorem~\ref{thm:PROMbound}}

\begin{proof}
The proof has two steps: an expression for $\hat\Sigma_e$ in terms of the discarded directions $Q$, and its factorization.  Both use the block form of $\Phi$ in the coordinates $U=[P\ \ Q]$,
\begin{equation*}
N=U^T\Phi U=\begin{bmatrix}sI_r-A_{PP}&-A_{PQ}\\-A_{QP}&sI-A_{QQ}\end{bmatrix},
\end{equation*}
with blocks $N_{ij}$.  Since $U$ is orthogonal, $N^{-1}=U^T\Phi^{-1}U$; in particular, $(N^{-1})_{22}=Q^T\Phi^{-1}Q$.

\emph{Error transfer matrix.}  The feedthrough $D$ cancels, and $\hat\Sigma_r=CP\,N_{11}^{-1}P^TB$, so $\hat\Sigma_e=C(\Phi^{-1}-PN_{11}^{-1}P^T)B$.  With $S=N_{22}-N_{21}N_{11}^{-1}N_{12}$, the partitioned inverse \cite{bernstein2009matrix} gives $N^{-1}-\operatorname{diag}(N_{11}^{-1},0)=[-N_{11}^{-1}N_{12};\,I]\,S^{-1}\,[-N_{21}N_{11}^{-1}\ \ I]$.  The two outer factors times $S^{-1}$ are the second block column and row of $N^{-1}$, and $(N^{-1})_{22}=S^{-1}$, so the right-hand side equals $(N^{-1})_{:,2}\,(N^{-1})_{22}^{-1}\,(N^{-1})_{2,:}$.  Multiplying by $U$ on the left and $U^T$ on the right, and using $UN^{-1}=\Phi^{-1}U$, gives
\begin{equation*}
\Phi^{-1}-PN_{11}^{-1}P^T=\Phi^{-1}Q(Q^T\Phi^{-1}Q)^{-1}Q^T\Phi^{-1},
\end{equation*}
which for symmetric positive-definite $\Phi$ is Khatri's lemma \cite{khatri1966note}.  Substituting into $\hat\Sigma_e$ proves \eqref{eq:PROM error TFM}.

\emph{Factorization.}  Let $\hat\Delta=Q^T\Phi^{-1}Q$, $\hat\Theta=C\Phi^{-1}Q\hat\Delta^{-1}$, and $\hat\Gamma=\hat\Delta^{-1}Q^T\Phi^{-1}B$, so that $\hat\Sigma_e=\hat\Theta\hat\Delta\hat\Gamma$ by \eqref{eq:PROM error TFM}; $\hat\Delta$ has the realization $(A,Q,Q^T,0)$ directly.  The off-diagonal blocks of $NN^{-1}=I$ and $N^{-1}N=I$ give $(N^{-1})_{12}=-N_{11}^{-1}N_{12}(N^{-1})_{22}$ and $(N^{-1})_{21}=-(N^{-1})_{22}N_{21}N_{11}^{-1}$.  Inserting $I=PP^T+QQ^T$ next to $\Phi^{-1}$, and using $P^T\Phi^{-1}Q=(N^{-1})_{12}$ and $Q^T\Phi^{-1}P=(N^{-1})_{21}$, yields
\begin{align*}
\Phi^{-1}Q\hat\Delta^{-1}&=P(sI_r-A_{PP})^{-1}A_{PQ}+Q,\\
\hat\Delta^{-1}Q^T\Phi^{-1}&=A_{QP}(sI_r-A_{PP})^{-1}P^T+Q^T.
\end{align*}
Hence $\hat\Theta=CP(sI_r-A_{PP})^{-1}A_{PQ}+CQ$ and $\hat\Gamma=A_{QP}(sI_r-A_{PP})^{-1}P^TB+Q^TB$, the transfer matrices of $\Theta$ and $\Gamma$ in \eqref{prod_form_elements}.
\end{proof}

\section*{Acknowledgment}
During the preparation of this work, the authors used Anthropic Claude to structure some of the proofs, to help develop the numerical case study, and to assist with editing the manuscript. After using this tool/service, the authors reviewed, restructured, and edited the content as needed and take full responsibility for the content of the publication.

\IEEEtriggeratref{40}

\end{document}